\documentclass[11pt,a4paper]{amsart}

\usepackage{ae} 
\usepackage[T1]{fontenc}
\usepackage[cp1250]{inputenc}
\usepackage{amsmath}
\usepackage{amssymb}
\usepackage{color}
\usepackage[colorlinks]{hyperref}
\usepackage[slovene, english]{babel}
\usepackage[dvips]{graphicx}
\usepackage{pstricks}

\def\proof{{\bf Proof. }}
\def\remark{{\hskip-4mm \bf Remark. }}

\newtheorem{dfn}{Definition}
\newtheorem{prop}[dfn]{Proposition}
\newtheorem{thm}[dfn]{Theorem}
\newtheorem{lemma}[dfn]{Lemma}

\newtheorem{cor}[dfn]{Corollary}
\newtheorem{corr}[dfn]{Corollary}

\renewcommand{\thedfn}%
{\arabic{dfn}}

\renewcommand{\rightsquigarrow}%
{ \hspace{4mm}\rput{0}(0,.1){\diagup \hskip-1.3mm \searrow}\hspace{4mm} }

\newcommand{\kdr}%
{\hspace{1mm}\rput{-45}(.05,1){\leadsto} \hspace{2mm}}
\newcommand{\nee}%
{\nearrow\hskip-4mm{{}^e}\hspace{3mm}}
\renewcommand{\see}%
{\searrow\hskip-3mm{{}^e}\hspace{2mm}}

\newcommand{\rel}%
{\emph{\hspace{0.1 cm}rel } }

\def\eps{\varepsilon}

\def\NN{\mathbb{N}}

\def\RR{\mathbb{R}}
\def\ZZ{\mathbb{Z}}

\def\i{^{-1}}
\newcommand{\rp}%
{ is relatively prime to }
\newcommand{\nrp}%
{ is not relatively prime to }

\def\di{\partial}                

\def\a{\alpha}
\def\b{\beta}

\def\d{\delta}

\def\f{\varphi}

\renewcommand{\theenumi}%
{\roman{enumi}}
\renewcommand{\labelenumi}%
{(\theenumi)}

\input pstricks.tex
\input xy
\xyoption{all}

\title{HOMOTOPICAL SMALLNESS AND CLOSENESS}
\author{ \v Ziga Virk }
\address{Faculty of Mathematics and Physics, University of Ljubljana, Jadranska 19, Ljubljana 1000, Slovenia
}
\email{virk@fmf.uni-lj.si}
\subjclass[2000]{54C20, 54D05, 54F15, 54G15, 54G20;  sets of homotopy classes of maps; small map; close map; homotopically Hausdorff space; generalized universal covering space; Spanier group;}
\thanks{The author would like to thank Matija Cencelj, Jerzy Dydak and Ale\v s Vavpeti\v c for their valuable comments.}

\begin{document}

\begin{abstract}
The aim of this paper is to introduce the concepts of homotopical smallness and closeness. These are the properties of homotopical classes of maps that are related to recent developments in homotopy theory and to the construction of universal covering spaces for non-semilocally simply connected spaces, in particular to the properties of being homotopically Hausdorff and homotopically path Hausdorff. The definitions of notions in question and their role in homotopy theory are supplemented by examples, extensional classifications, universal constructions and known applications.
\end{abstract}

\maketitle

\section{Introduction}
The concepts of homotopical smallness and closeness are related to various versions of the property of being homotopically Hausdorff, which have been introduced and studied in \cite{CC1}, \cite{Re} and \cite{Z}.

\begin{dfn} \label{homotopHauss}
A space $X$ is called:
\begin{enumerate}
\item  (weakly)  \textbf{homotopically Hausdorff} if for every $x_0 \in X$ and for every  non-trivial $\gamma \in \pi_1(X,x_0)$ there exists a neighbourhood $U$ of $x_0$ such that no loop in $U$ is homotopic (in $X$) to $\gamma$ rel.\ $x_0$.
\item  \textbf{strongly homotopically Hausdorff} if for every
   $x_0 \in X$ and for every essential closed curve $\gamma \in X$ there is a neighbourhood of $x_0$ that contains no closed curve freely homotopic (in $X$) to $\gamma$.
\item \textbf{homotopically path-Hausdorff}
      if for every path $w \colon [0,1] \to X$ with $w(0)=P$ and $w(1)=Q$ and every non-trivial
      homotopy class $\alpha \in \pi_1(X,P)$ there exist finitely many open sets
      $U(P_1),\ldots,U(P_k), \quad (P_1=P$ and $P_k=Q)$ covering $w([0,1])$
      such that for a suitable partition $0=t_0 < t_1 < t_2 < t_3 <\ldots< t_k=1,\quad
      U(P_j)$ covers $w([t_{j-1},t_j])$, $P_j\in w([t_{j-1},t_j])$ and such that for any  path $v\colon [0,1] \to X$ that satisfies
      $v(0)=Q $, $ v(1)=P $, $ P_{k-j+1} \in v([t_{j-1},t_j])$ and
      $v([t_ {j-1} , t_j]) \subset U(P_ {k-j+1})~\forall j$
      the concatenation
      of $w$ and $v$ does not belong to the homotopy class $\alpha$.
\end{enumerate}
\end{dfn}

Note that the property of being homotopically Hausdorff is weaker than both the property of being strongly homotopically Hausdorff and the property of being homotopically path Hausdorff. These are separation properties for homotopical classes of maps and play a significant role in homotopy theory for locally wild spaces, for example, spaces which are not semilocally simply connected, etc.  Good examples of such spaces are Hawaiian earring (denoted by $HE$) and Harmonic archipelago (denoted by $HA$).

The $HE$ is a countable metric wedge of circles (circular loops) whose diameters tend to zero, i.e.,
$$
HE:=\bigcup_{i\in \ZZ_+} S^1((\frac{1}{i},0),\frac{1}{i})\subset \RR^2
$$
where $S^1(C,r)\cong S^1$ is the circle in $\RR^2$ with center $C$ and radius $r$. Circles $S^1((\frac{1}{i},0),\frac{1}{i})$ are equipped with the positive (respectively negative) orientation and denoted by $l_i$ (respectively $l_i^-$). The intersection of all circles is denoted by $0$. It turns out that $HE$ is homotopically Hausdorff but not semilocally simply connected.

\begin{figure}
\begin{center}
\includegraphics{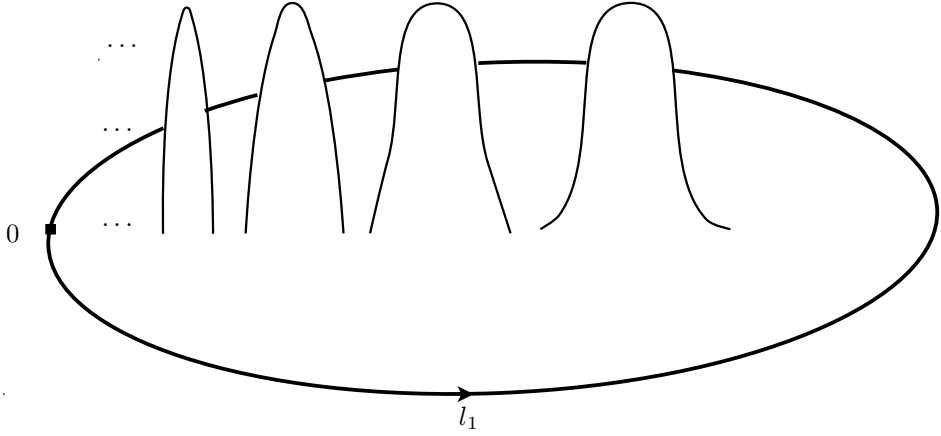}
\end{center}
\caption{Harmonic archipelago.}
\label{HA}
\end{figure}

The Harmonic archipelago was defined in \cite{BS} and studied in \cite{Fa}. In order to construct it begin with $HE\subset \RR^2\times \{0\}\subset \RR^3$. For each pair of consecutive loops $(l_i, l_{i+1})$ attach the disc $B_i^2$ in the following way: identify the boundary $\di B^2_i$ with the loop $l_i * l_{i+1}^-$ and stretch the interior of $B^2_i$ up so that one of its interior points (called the peak point of $B^2_i$) is at height $1$. The situation is presented in Figure \ref{HA} where discs $B^2_i$ are represented by "bumps". It is easy to see that $HA$ is not homotopically Hausdorff.

A detailed study of relationship between properties of Definition \ref{homotopHauss} is presented in \cite{Re} and \cite{4aut}. The distinction between them is demonstrated by spaces $Y,Y',Z,Z'$ of \cite{4aut}. Essential parts of these spaces turn out to be a generic spaces where certain properties of smallness and closeness occur.

Properties of Definition \ref{homotopHauss} arise in connection to the universal path space.
\begin{dfn}\label{ups}
    Let $(X,x_0)$ be a pointed path connected space. The \textbf{universal path space} $\widehat X$ is the set of equivalence classes of paths $\a\colon [0,1]\to X, \a(0)=x_0$ under the following equivalence relation: $\a \sim \b$ iff $\a(1)=\b(1)$ and the concatenation $\a*\b^-$ (where $\b^-(t):=\b(1-t)$) is homotopic to  a constant path at $x_0$, denoted by $1_{x_0}$. The space $\widehat X$ is given a topology generated by the  sets
    $$
    N(U,\a):=\{\b \mid \b \simeq \a * \eps, \eps\colon ([0,1],0)\to (U, \a(1))\}
    $$
    where $U$ is an open neighborhood of $\a(1)\in X$. The natural endpoint projection $\hat p\colon \widehat X\to X$  is called the \textbf{endpoint map}.

    Universal path space is called a \textbf{universal covering space} if the endpoint projection has the unique path lifting property.
\end{dfn}

The following are well known facts that appear in \cite{CC1}, \cite{Re}, \cite{4aut}, \cite{Vi1} and \cite{Z}.

\begin{prop}
    Let $(X,x_0)$ be a path connected space.
\begin{enumerate}
  \item $X$ is semilocally simply connected iff the fibers $\hat p^{-1}(x)\subset \widehat X$ of the endpoint projection are discrete subspaces for all $x\in X$.
  \item $X$ is homotopically Hausdorff iff the fibers $\hat p^{-1}(x)\subset \widehat X$ of the endpoint projection are Hausdorff subspaces for all $x\in X$.
  \item If $\widehat X$ is a universal covering space then $X$ is homotopically path Hausdorff.
  \item If $X$ is homotopically Hausdorff and $\pi_1(X,x_0)$ is countable then $\widehat X$ is a universal covering space.
  \item If $X$ is homotopically path Hausdorff then $\widehat X$ is a universal covering space.
\end{enumerate}
\end{prop}

The property of being homotopically Hausdorff is closely related to small loops, which were introduced and studied in \cite{Vi1}. In this paper we extend the approach of \cite{Vi1} in order to define smallness and closeness for a wider class of maps and relate new concepts to existing examples and properties. As a result we obtain the following classification.

\begin{thm} \label{classifOFhomotopt2}Let $X$ be a path connected space.
\begin{enumerate}
  \item A space $X$ is homotopically Hausdorff if it contains no non-trivial pointed small loop.
  \item A space $X$ is strongly homotopically Hausdorff if it contains no non-trivial free small loop.
  \item A locally path connected space $X$ is homotopically path-Hausdorff if there is no pair of paths in $X$ that are close relatively to the endpoints of the interval.
\end{enumerate}
\end{thm}

Statements (i) and (ii) are apparent from Definitions \ref{smallMapPointed} and \ref{SmallMap}. Statement (iii) is the content of Proposition \ref{ClosePathsAndHomotopPathHaus}.

\section{Technical preliminaries}

We introduce several notions that will be used in the course of the paper. The following definition is a generalization of a classical concept of an absolute extensor which will be used to classify certain cases of smallness and closeness.

\begin{dfn}
Let $A\subseteq X$ be a closed subspace and let $Y$ be any topological space. Space $Y$ is an \textbf{absolute extensor} for the inclusion  $A\hookrightarrow X$ $[$notation: $(A\hookrightarrow X)\tau Y$ or $Y \in AE(A\hookrightarrow X)]$ if every map $A\to Y$ extends over $X$.
\end{dfn}

Note that a path connected space $Y$ is simply connected iff it is an absolute extensor for the inclusion $\di B^2 \hookrightarrow B^2$.

The notion of an m-stratified space as defined in \cite{Vi1} is a description of construction rather than the property of a space as every space is m-stratified. It mimics the structure of $CW$-complexes by building spaces through attachment of smaller pieces via quotient maps.

\begin{dfn}\label{mStratified}
Let $\{Y_i, A_i\}_{i\geq 0}$  be a countable collection of pairs of spaces where $A_i\subseteq Y_i$ is closed for every $i.$ Topological space $X$ is an \textbf{m-stratified} (map stratified) space with parameters $\{Y_i,A_i\}_i$ if  it is homeomorphic to the direct limit of spaces $\{X_i\}_{i\geq 0}$ where spaces $X_i$ are defined inductively as
\begin{itemize}
  \item $X_0:=Y_0$,
  \item $X_i:= X_{i-1}\cup_{f_i}Y_i$ for some maps $f_i\colon A_i \to X_{i-1}$.
\end{itemize}
The sets $Y_i$ are called \textbf{m-strata}.
\end{dfn}

When applying the construction of an m-stratification we will usually adopt the notation of Definition \ref{mStratified}. Lemma \ref{CompactInQuotient} has origins in the theory of $CW$ complexes presented in \cite{Ha}. It describes the behavior of compact subsets with respect to an m-stratification.

\begin{lemma} \cite{Vi1}
\label{CompactInQuotient}
Suppose  $Y$ is an m-stratified space so that m-strata $Y_i$ can be decomposed as $Y_i=\coprod_j Y_i^j$ where $Y_i^j\subset Y_i$ are open regular subspaces (i.e. open subspaces which are regular topological spaces). Let $K\subset Y$ be a  compact space. Define $X_i^j$ to be the image of $Y_i^j$ in $Y$. Then $K$ is contained in a finite union of subsets $X_i^j\subset Y$.
\end{lemma}

Another important property is related to extensions of maps. Any synchronized collection of maps on m-stratas induces a continuous map on $Y$.

\begin{lemma}\label{extStrat}
Let $Y$ be an m-stratified space and let $g_i\colon Y_i \to Z$ be a collection of maps satisfying $g_i|_{A_i}=g_i|_{f_i(A_i)}\circ f_i$. Then maps $g_i$ induce a continuous map on $Y$.
\end{lemma}

The notion of a universal Peano space was defined in \cite{D2}. It allows us to study certain properties of a non-locally path connected space.

\begin{dfn} \label{PeanoSpace}
Let $X$ be a path connected space. The \textbf{universal Peano space} (or peanification) $PX$ of $X$ is the set $X$ equipped with a new topology, generated by all path components of all open subsets of the existing topology on $X$. The \textbf{universal Peano map} is the natural bijection $p\colon PX\to X$.
\end{dfn}

Note that $PX$ is locally path connected. As an example, the peanification of the Warsaw circle is a semi-open interval. The name "universal Peano map" refers to the universal map lifting property for locally path connected spaces.

\begin{prop} \cite{D2} \label{ZZZ}
Let $Y$ be a locally path connected space. Then every map $f\colon Y\to X$ uniquely lifts to a map $f'\colon Y\to PX$.
\end{prop}
$$
\xymatrix{&PX \ar[d]^p\\
Y \ar@{-->}[ur]^{f'} \ar[r]_f&X
}
$$

\proof Since $p$ is bijection the only possible choice for $f'$ is $p\- f$. Let us prove it is continuous. Choose $y\in Y$ and let $x=f(y)$. Every open neighborhood $U'\subset PX$ of $x'=p\-(x)\in PX$ is a path component of an open neighborhood $U\subset X$ of $x\in X$. The preimage $f\-(U)$ is an open neighborhood of $y$ which contains an open path connected neighborhood $W$ of $y$ as $Y$ is locally path connected. Then $f(W)\subset U$ is path  connected and contains $x$ hence $f'(W)$ is contained in $U'$. \hfill$\blacksquare$
\bigskip

If $Y$ is locally path connected then so is $Y\times [0,1]$ which yields  the following corollary.

\begin{cor} \label{MapsToPeano} Let $Y$ be  a locally path connected space and let $X$ be a path connected space.
\begin{enumerate}
  \item The set of homotopy classes of maps $[Y,X]$ is in natural bijection with $[Y,PX].$
  \item  The set of homotopy classes of maps $[Y,X]_\bullet$ in the pointed category is in natural bijection with $[Y,PX]_\bullet.$
  \item $\pi_k(X)=\pi_k(PX)$, for all $k\in \ZZ^+$.
  \item $H_k(X)=H_k(PX)$, for all $k\in \ZZ^+$.
\end{enumerate}
\end{cor}

Corollary \ref{ZZZ} implies that paths and homotopies between paths, on the base of which the universal path space is defined, are the same in $X$ and $PX$.

\begin{cor}\label{ZZZZ}
The universal path spaces $\widehat X$ and $\widehat {PX}$ are homeomorphic for every path connected space $(X,x_0)$.
\end{cor}

Given a path connected space $X$ which is not locally path connected, Corollary \ref{ZZZZ} describes the information about the space which is lost in the construction of the universal path space. Essentially it is the same information that is lost in the construction of the universal Peano space. In particular, the universal Peano space $PX$ has the same homotopy and homology groups as $X$ but may have different shape group. On the other hand, spaces $\widehat X$ and $PX$ are locally path connected even if the space $X$ is not.

\section{Homotopical smallness}

The definition of small maps first appeared in \cite{Vi1} in the form of small loops.

\begin{dfn} \label{slg}
A loop $\a\colon (S^1,0)\to (X,x_0)$ is \textbf{small} iff there exists a representative of the  homotopy class $[\a]_{x_0}\in \pi_1(X,x_0)$ in every open neighborhood $U$ of $x_0$. A small loop is a \textbf{non-trivial small loop} if it is not homotopically trivial.
\end{dfn}

Griffiths' space of \cite{G} and  $HA$ of \cite{BS} are well known spaces with non-trivial small loops. Another example is the strong Harmonic Archipelago $SHA$. The topology of $SHA$ can be described in terms of m-stratified spaces with the following parameters (using the notation of Definition \ref{mStratified}):
$$
Y_0=HE, \quad, Y_i=B^2_i, \quad A_i=\di B^2_i=S^1_i, \quad f_i=l_i l^-_{i+1}\colon S^1_i\to HE.
$$

Both $HA$ and $SHA$ are obtained from $HE$ by attaching discs $B^2_i$ along loops $l_i l^-_{i+1}$. The difference is that in the case of $SHA$ an infinite collection of discs $\{B^2_i\}$ is attached to $HE$ by the quotient map (making it more natural as suggested by the proof of Proposition \ref{HHHH}), while in the case of $HA$ attachment is carried on in $\RR^3$ so that the resulting space $HA$ is metric. $SHA$ is a generic example of a non-trivial small loop in a first countable in the sense of the following proposition, which can be proved using Lemma \ref{extStrat}. Its generalization will be proven later.

\begin{prop} \label{HHHH}
Assume that   $x_0\in X$  has a countable basis of neighborhoods. A loop $\a\colon (S^1,0)\to (X,x_0)$ is small iff it extends to $F\colon (SHA,0)\to (X,x_0)$ where $l_1\colon (S^1,0)\hookrightarrow (SHA,0)$ is the boundary loop.
\end{prop}
$$
\xymatrix{(S^1,0) \ar@{^(->}[dr] \ar[r]^\a& (X,x_0)\\
&(SHA,0) \ar@{-->}[u]
}
$$

The same proposition can be proven for $HA$ instead of $SHA$ as well but the proof is somewhat more complicated. Another construction related to small loops are small loop spaces as defined, constructed and studied in \cite{Vi1}.

\begin{dfn}
A non-simply connected space $X$ is a \textbf{small loop space} if for every $x\in X$,
every loop $\a\colon (S^1,0)\to (X,x)$ is small.
\end{dfn}

The following subsections generalize the notion of smallness and accompanying constructions to a general case in various categories.

\subsection{Homotopical smallness in unpointed category}
This subsection is devoted to homotopical smallness of arbitrary spaces in unpointed category. All homotopies and maps are considered to be unpointed (i.e. spaces have no basepoint and homotopies need not preserve any point).

We start with a definition of smallness in the unpointed category. The absence of a basepoint implies that we should specify a point at which we would like to consider smallness. By smallness we mean the property of being able to find a homotopic representative of a map in every neighborhood of a point.

\begin{dfn}\label{SmallMap}
A map $f\colon Y\to X$  is (homotopically freely) \textbf{small} at $x\in X$  (in unpointed topological category) if for each open neighborhood $U$ of $x$ there is a (free) homotopy $H\colon Y \times [0,1]\to X$  so that $H|_{Y \times \{0\}}=f$ and $H|_{Y \times \{1\}}(Y)\subset U$. A small map is a \textbf{non-trivial small map} if it is not homotopically trivial.
\end{dfn}

\begin{prop}
Suppose $f\colon Y\to X$ is a  small map  at $x\in X$ and $g\colon f(Y)\to Z$ is a map. If $g$ extends over $X$ then $gf\colon Y\to Z$ is a small map at $g(x)$.
\end{prop}

For the rest of this chapter we will assume $S$ to be a directed set with no maximal element (hence $S$ is infinite) and the smallest (initial) element $s_0$, unless otherwise stated.  Definition \ref{FSOS} introduces a generic examples of small maps which classify all small maps in terms of extension theory.

\begin{dfn} \label{FSOS}
The \textbf{Sydney opera space} of $S$ with respect to the space $Y$ (in the topological category) [notation: $FSO_Y(S)$] is a space constructed in the following way.

Take a disjoint union $\coprod_{s\in S} Y_s$ of copies of space $Y$, one copy for each element in $S$. Upon this union attach spaces $W_s:=Y\times [0,1]$  for each $s\in S\backslash\{s_0\}$, so that $Y\times\{0\}\subset W_s$ is identified with $Y_{s_0}$  and $Y\times\{1\}\subset W_s$ is identified with $Y_s$. Add another point $\{0\}$ to obtain the space $FSO_Y(S):=\bigcup_{s\in S \backslash \{s_0\}} W_s \cup \{0\}$ and define the following topology. The subset $U\subset FSO_Y(S)$ is open if either of the following is true:
\begin{enumerate}
  \item $0\notin U$ and $U$ is open in  $W_s, \forall s\in S \backslash \{s_0\}$,
  \item $0\in U$, $U$ is open in  $W_s, \forall s\in S\backslash \{s_0\}$ and there exists $t_0\in S $ such that $Y_t\subset U, \forall t\geq t_0$.
\end{enumerate}
\end{dfn}

For a fixed directed set $S$ with the initial element $s_0$ the rule $Y\mapsto FSO_Y(S)$ is a functor on the category of the topological spaces. Space $FSO_Y(S)$ can be given various structures of an m-stratified space. The simplest one would start with $\{0\}\cup \coprod_{s\in S} Y_s$ (with appropriate topology as described in Definition \ref{FSOS}) upon which we attach  homotopies $W_s.$ Using the notation of Definition \ref{mStratified} the topology of $FSO_Y(S)$ can be expressed by the following parameters:
$Y_0= \{0\}\cup\coprod_{s\in S} Y_s$ (with topology described in Definition \ref{FSOS}),
$$
Y_1=\coprod_{s\in S\backslash\{s_0\}}(Y\times [0,1])_s, \quad A_1 = \coprod_{s\in S\backslash\{s_0\}}(Y \times\{0,1\})_s,
$$
$$
f_1|_{(Y \times\{0\})_s}=1_{ Y_0}, \quad f_1|_{(Y \times\{1\})_s}1_{Y_s}.
$$

Note that $0\in FSO_Y(S)$ is not path connected to $Y_{s_0}$.

\begin{lemma}
The natural inclusion $Y\to Y_{s_0}\subset FSO_Y(S)$ is small at $0$.
\end{lemma}

\proof Using homotopies $W_s$ we can homotope the inclusion into arbitrary neighborhood of $0$. \hfill $\blacksquare$
\bigskip

If $Y$ is contractible then the inclusion $Y\to Y_{s_0}\subset FSO_Y(S)$ is homotopically trivial. A necessary condition for such inclusion to be homotopically non-trivial is homotopical non-triviality of $Y$. Sufficient condition is given by Corollary \ref{endvatri}.

\begin{lemma}
\label{FactorFreeHomotopy}
Let $f\colon K\to FSO_Y(S)$ be a map from a compact space $K$ to a regular space $Y$. Then $f(K)$ is contained in the subspace
$$
\bigcup_{s\in T}W_s \cup \bigcup_{s\in S}Y_s \cup \{0\}
$$
where $T\subseteq S$ is some finite subset. Furthermore, such $f$ factors over $\cup_{s\in S}Y_s\cup \{0\}\hookrightarrow FSO_Y(S)$ up to homotopy.
\end{lemma}

\proof
The first part follows by Lemma \ref{CompactInQuotient}. To prove the second part consider a strong deformation retraction
$$
\bigcup_{s\in T}W_s \cup \bigcup_{s\in S}Y_s \cup \{0\} \to \bigcup_{s\in S-
T}Y_s\cup Y_{s_0}\cup\{0\}.
$$
\hfill $\blacksquare$

\begin{lemma}\label{FactorFreeHomotopy1}
Let $f\colon K\to Y_{s_0}\subset FSO_Y(S)$ be a map from a compact space $K$ to a regular space $Y$  and suppose $H\colon K\times [0,1] \to FSO_Y(S)$ is a homotopy so that $H|_{K\times \{0\}}=f$. Then $H(K\times [0,1])$ is contained in the subspace
$$
\bigcup_{s\in T}W_s
$$
where $T\subseteq S$ is some finite subset.
\end{lemma}

\proof
By Lemma \ref{FactorFreeHomotopy} the compact set  $H(K\times [0,1])$ is contained in
$$
A:=\bigcup_{s\in T}W_s \cup \bigcup_{s\in S}Y_s \cup \{0\}
$$
for some finite $T\subseteq S$. Note that $H(K\times [0,1])$ is connected to $Y_{s_0}$ by paths as  $H|_{K\times \{0\}} \subset Y_{s_0}$. Sets $\{0\}$ and $Y_t\cap (K\times [0,1])$ for $t\in S-T$ are not path connected to $Y_{s_0}$  in $A$ hence  $H(K\times [0,1])\subseteq \cup_{s\in T}W_s$. \hfill $\blacksquare$

\begin{cor}\label{endvatri}
Let $Y$ be a compact Hausdorff space which is not homotopically trivial. Then the natural inclusion $i\colon Y\to Y_{s_0}\subset FSO_Y(S)$ is a homotopically non-trivial small map at $0$.
\end{cor}

\proof
Suppose there is a homotopy $H$ in $FSO_Y(S)$ between the inclusion $i$ and a constant map. Because $Y$ is compact such homotopy is a compact map therefore its image is contained in $\cup_{t\in T}W_s$ where $T\subseteq S$ is a finite subset.
Space  $\cup_{t\in T}W_t$ can be naturally retracted to $Y_{s_0}$. Composing homotopy $H$ with such retraction we  contradict the fact that $Y$ is homotopically non-trivial.  \hfill $\blacksquare$
\bigskip

Similarly as in the case of small loops we can classify small maps in terms of extension theory.

\begin{prop} Let $S$ be a directed set with the smallest element $s_0$ so that for a point $x\in X$ there is a basis $\{U_s\}_{s\in S}$ that satisfies  $(U_s\subseteq U_t) $ iff $(s\geq t)$. Map $f\colon Y_{s_0}\to X$ is small at $x\in X$ iff it extends over $FSO_Y(S)$ to a map $F$ so that $F(0)= x$.
\end{prop}

\proof
We only have to prove one direction. Suppose the map $f\colon Y_{s_0}\to X$ is  small at $x\in X$. For each $s\in S$ there is a homotopy $H$ between $f$ and a map whose image is contained in $U_s$. Use such homotopy to naturally define the map $F$  on $W_s$ and additionally define $F(0):=x$. This rule defines a continuous map on $FSO_Y(S)-\{0\}$  as the topology on it  is quotient. The preimage $F^{-1}(U_t)$ of any basic open neighborhood $U_t$ of $x$ is open in  $W_s, \forall s\in S,$ and contains entire $Y_s$ for all $s\geq t$ therefore it is open in $FSO_Y(S)$. Hence the extension $F$ is continuous.
\hfill $\blacksquare$

\begin{dfn} \label{SmallKSpace}
Let $Y$ be a topological space. Space $X$ is called a \textbf{small $Y-$space} if the following conditions hold:
\begin{enumerate}
  \item There exists a map $Y\to X$ which is not homotopically trivial.
  \item Every map $Y\to X$ is small at every  $x\in X$.
\end{enumerate}
\end{dfn}

\begin{cor}
Let $S$ be a directed set with the smallest element $s_0$ so that for $x\in X$ there is a basis $\{U_s\}_{s\in S}$ that satisfies  $(U_s\subseteq U_t) $ iff $(s\geq t)$.
Suppose  there exists a map $Y\to X$ which is not homotopically trivial. Space $X$ is a small $Y-$space iff $\Big((Y_{s_0}\cup \{0\})\hookrightarrow FSO_Y(S)\Big)\tau X. $
\end{cor}

The rest of this subsection is devoted to the existence of a small $Y-$space for  non-contractible compact Hausdorff spaces. With these properties we can imitate the construction of a small loop space of \cite{Vi1}. The Hausdorff property of $Y$ implies that $Y$ and $FSO_Y(\NN)$ are regular spaces which allows Lemma \ref{CompactInQuotient} to be used in the case of $FS(Y)$. The following definition introduces a generic example of a small $Y$-space.

\begin{dfn}\index{FS(Y) space}
Let $Y$ be a topological space. The space $\textbf{FS(Y)}$ is an m-stratified space with
$$
Y_0=FSO_Y(\NN), \qquad S_i:=\{(g,x);\quad  g\colon Y\to X_{i-1}, \quad x\in X_{i-1}\},
$$
$$
Y_i:= \coprod_{(g,x)\in S_i} FSO_Y(\NN)_{g,x}, \qquad A_i:= \coprod_{(g,x)\in S_i} (Y_0\cup \{0\})_{g,x}
$$
$$
f_i({0}_{g,x})=x, \qquad f_i|_{(Y_0)_{g,x}}=g,
$$
where $(Y_0)_{g,x}\subset FSO_Y(\NN)_{g,x}$ is the initial copy of $Y$ in $FSO_Y(\NN)_{g,x}$.
\end{dfn}

\begin{lemma}
Suppose $Y$ is a compact Hausdorff space. Every map $f\colon Y \to FS(Y)$  is small.
\end{lemma}

\proof Choose any $x\in FS(Y)$. By Lemma \ref{CompactInQuotient} there is $k\in \NN$ so that $f(Y)\subset X_k$ and $x\in X_k$ according to m-stratification of $FS(Y)$. Then $f$ can be made small in $X_{k+1}$ via the attached space $FSO_Y(\NN)_{f,x}$. \hfill $\blacksquare$

\begin{prop}
Suppose $Y$ is a compact Hausdorff space which is not homotopically trivial. Then the natural  inclusion $f\colon Y\to Y_{s_0}\subset X_0 \subset FS(Y)$ is homotopically non-trivial
\end{prop}

\proof Suppose there is a homotopy $H$ taking $f$ to a constant map. Using Lemmas \ref{CompactInQuotient} and \ref{FactorFreeHomotopy1}  one can construct a retraction of $H(Y\times I)$ to $Y_{s_0}$. Composing such retraction with  $H$ would imply that $Y$ is homotopically trivial, which is a contradiction. \hfill $\blacksquare$

\begin{corr}
If the space $Y$ is compact Hausdorff and homotopically non-trivial then $FS(Y)$ is a small $Y-$space.
\end{corr}

Space $FS(Y)$ is a universal example of a small $Y$-space in the following way.

\begin{prop}
Suppose $f\colon Y_{s_0}\to X$ is a map to a small $Y-$ space $X$ where $Y_{s_0}\subset X_0 \subset FS(Y)$. Then $f$ extends over $FS(Y).$
\end{prop}

\proof Follows from Lemma \ref{extStrat}. \hfill$\blacksquare$
\bigskip

\subsection{Homotopical smallness in pointed category}

The aim of this subsection is to develop similar results for smallness in the pointed category. Recall that pointed homotopy is a homotopy that fixes the base point. All spaces, maps and homotopies of this section are considered to be in the pointed topological category.

\begin{dfn}\label{smallMapPointed}
A map $f\colon (Y,y_0)\to (X,x_0)$ between pointed topological spaces is (homotopically) \textbf{ small}  (in the pointed topological category) if for each open neighborhood $U$ of $x_0$ there exists a (pointed) homotopy (i.e. $H\colon (Y\times[0,1],(y_0,0))\to (X,x_0), \quad H(y_0,t)=x_0, \quad \forall t\in [0,1]$) so that  $H|_{Y\times \{0\}}=f$ and  $H|_{Y\times \{1\}}(Y\times \{1\})\subset U$. A small map is a \textbf{non-trivial small map} if it is not homotopically trivial.
\end{dfn}

\begin{prop} Let  $f\colon (Y,y_0)\to (X,x_0)$ be a  small map and let  $g\colon (f(Y),x_0)\to (Z,z_0)$ be a map. If $g$ extends over $X$ then $gf\colon (Y,y_0)\to (Z,z_0)$ is small.
\end{prop}

The following definition introduces an analogue of the $FSO$ spaces in the pointed category. Recall that $S$ is assumed to be a directed set with no maximal element and the smallest element $s_0$.

\begin{dfn} \index{Sydney opera space} \index{SO${_Y}$(S) space}\label{SOS}
The \textbf{ Sydney opera space} [notation: $SO_Y(S)$] of $S$ with respect to the space $(Y,y_0)$ (in the pointed topological category) is a topological space constructed in the following way.

Consider the wedge $\vee_{s\in S} (Y_s,y_0)$ of $|S|$ copies of the space $Y$, one copy for each element of $S$, obtained by identifying the base points $y_0$ of the spaces $(Y_s,y_0)$. Define its basepoint to be the wedge point and denote it by $y_0$ as well. On this wedge attach the spaces $W_s:=Y\times [0,1]$  for each $s\in S\backslash \{s_0\}$, so that $Y\times\{0\}\subset W_s$ is identified with $Y_{s_0}$,   $Y\times\{1\}\subset W_s$ is identified with $Y_s$ and $\{y_0\}\times I$ is identified with $y_0\in \vee_{s\in S\backslash \{s_0\}} (Y_s,y_0)$. Define $SO_Y(S)$ to be the set $\cup_s W_s$ with the following topology. A subset $U\subset SO_Y(S)$ is open if either of the following is true:
\begin{enumerate}
  \item $y_0\notin U$ and $U$ is open in  $W_s, \forall s\in S\backslash \{s_0\}$,
  \item $y_0\in U$, $U$ is open in  $W_s, \forall s\in S\backslash \{s_0\}$ and there is $t_0\in S$ so that $Y_t\subset U, \forall t\geq t_0$.
\end{enumerate}
\end{dfn}

For a fixed directed set $S$ with the initial element $s_0$ the rule $Y\mapsto SO_Y(S)$ is a functor in the category of the pointed topological spaces. The space $SO_Y(S)$ is a natural quotient of $FSO_Y(S)$ and can be given various structures of an m-stratified space. The simplest one is almost identical to the one of $FSO_Y(S)$. Start with the wedge $\vee_{s\in S} (Y_s,y_0)$ (with appropriate topology as described in Definition \ref{SOS}) upon which we attach homotopies $W_s.$ Using the notation of Definition \ref{mStratified} the topology of $SO_Y(S)$ can be expressed by the following parameters:
$Y_0= \coprod_{s\in S} Y_s$ (with topology described in Definition \ref{SOS}),
$$
Y_1=\coprod_{s\in S\backslash\{s_0\}}(Y\times [0,1])_s, \quad A_1 = \coprod_{s\in S\backslash\{s_0\}}(Y \times\{0,1\}\cup \{y_0\}\times [0,1])_s,
$$
$$
f_1|_{(Y\times \{0\})_s}=1_{Y_0}, \quad f_1|_{(Y\times \{1\})_s}=1_{Y_s}\quad f_1\big((\{y_0\}\times [0,1])_s\big)= \{y_0\}.
$$

The topology of an m-stratified space implies that the natural inclusion $(Y,y_0)\cong (Y_{s_0},y_0)\subset (SO_Y(S),y_0)$ is small. The nature of compact subsets implies that such inclusion homotopically non-trivial if $Y$ is not contractible.

\begin{lemma}
Let  $S$ be a directed set with the smallest element $s_0$. The natural inclusion $(Y,y_0)\to (Y_{s_0},y_0)\subset (SO_Y(S),y_0)$ is small.
\end{lemma}

\proof Use homotopies $W_s$. \hfill $\blacksquare$

\begin{lemma}
\label{FactorHomotopy}
Suppose $f\colon (K,k_0)\to (SO_Y(S),y_0)$ is a map defined on a compact Hausdorff space $(K,k_0)$. Then $f(K)$ is contained in the subspace
$$
\bigcup_{s\in T}W_s \cup \bigcup_{s\in S}Y_s
$$
where $T\subset S$ is a finite subset. Furthermore, such $f$ factors over $(\cup_{s\in S}Y_s,y_0)\hookrightarrow (SO_Y(S),y_0)$ up to homotopy.
\end{lemma}

\proof
The first part follows by Lemma \ref{CompactInQuotient}. To prove the second part consider the strong deformation retraction
$$
\bigcup_{s\in T}W_s \cup \bigcup_{s\in S}Y_s  \to \bigcup_{s\in S-
T}Y_s\cup Y_{s_0}.
$$
\hfill $\blacksquare$

\begin{cor}
Let $(Y,y_0)$ be a compact Hausdorff, homotopically non-trivial space.  The natural inclusion $i\colon (Y,y_0)\to (Y_{s_0},y_0)\subset (SO_Y(S),y_0)$ is a homotopically non-trivial small map.
\end{cor}

\proof Suppose there is a pointed homotopy in $SO_Y(S)$ between $i$ and a constant map. The space $(Y,y_0)$ is compact hence the image of such homotopy is contained in $\bigcup_{s\in T}W_s \cup \bigcup_{s\in S-T}Y_s$ where $T\subseteq S$ is some finite subset. The subspace  $\cup_{t\in T}W_T\subset SO_Y(S)$ can be retracted to $Y_{s_0}$ and the subspace $\bigcup_{s\in S-T}Y_s\subset SO_Y(S)$  can be retracted to $y_0$. Composing the homotopy with these retractions we obtain a contraction of $(Y,y_0)$, a contradiction.  \hfill $\blacksquare$

\begin{prop}
Let $S$ be a directed set with the smallest element $s_0$ so that $x\in X$ has a basis of neighborhoods $\{U_s\}_{s\in S}$ that satisfies  $(U_s\subseteq U_t) $ iff $(s\geq t)$. The map $f\colon (Y_{s_0},y_0)\to (X,x_0)$ is small iff it extends over $SO_Y(S)$.
\end{prop}

\proof We only have to prove one direction. Suppose the map $f\colon (Y_{s_0},y_0)\to (X,x_0)$ is small. For each $s\in S$ there is a  homotopy $H$ between $f$ and a map with its image contained in $U_s$. Use such homotopy to naturally define a map on $W_s$. With this rule we have defined a continuous map on $SO_Y(S)-\{y_0\}$ by the definition of topology. The preimage of any basic open neighborhood $U_t ,t\in S$ of $x$ is open in  $W_s, \forall s\in S,$ and contains all $W_s, s\geq t$ therefore it is open. Hence the extension is continuous. \hfill $\blacksquare$
\bigskip

Definition \ref{SmallKSpacePointed} introduces a small $Y-$space. It is followed by the extensional classification and a construction of such space. Small $Y-$space is a generalization of a small loop space which was shown to have interesting properties in \cite{Vi1}.

\begin{dfn}  \label{SmallKSpacePointed}
Let $(Y,y_0)$ be a topological space. We call $X$ a \textbf{small $Y-$space} if the following conditions hold:
\begin{enumerate}
  \item There exists a map $f\colon (Y,y_0)\to (X,f(y_0))$ which is not homotopically trivial.
  \item Every map $f\colon (Y,y_0)\to (X,f(y_0))$ is small.
\end{enumerate}
\end{dfn}

\begin{corr}
Let $X$ be a topological space and suppose $S$ is a directed set with the smallest element $s_0$ so that for every $x\in X$ there exists a basis $\{U_s\}_{s\in S}$ of open neighborhoods of $x$ satisfying  $(U_s\subseteq U_t) $ iff $(s\geq t)$. Suppose there exists a homotopically non-trivial map $f\colon (Y,y_0)\to (X,f(y_0))$. Space $X$ is a small $Y-$space iff $\Big(Y_{s_0}\hookrightarrow SO_Y(S)\Big)\tau X. $
\end{corr}

\begin{dfn}  \index{S(Y) space}
Let $(Y,y_0)$ be a topological space. The space $\textbf{S(Y)}$ is an m-stratified space with
$$
Y_0=SO_Y(\NN), \qquad S_i:=\{g; \quad g\colon Y\to X_{i-1}\},
$$
$$
Y_i:= \coprod_{g\in S_i} SO_Y(\NN)_{g}, \qquad A_i:= \coprod_{g\in S_i} (Y_0)_{g}, \qquad f_i=\coprod_{g\in S_i}g,
$$
where $(Y_0)_{g,x}\subset SO_Y(\NN)_{g,x}$ is the initial copy of $Y$ in $SO_Y(\NN)_{g}$.
\end{dfn}

\begin{lemma}
If $(Y,y_0)$ is a compact Hausdorff space then every  $f\colon (Y,y_0) \to (S(Y),f(y_0))$ is small.
\end{lemma}

\proof
By Lemma \ref{CompactInQuotient} there exists $k\in \NN$ so that $f(Y)\subset X_k$. Pointed map $f\colon (Y,y_0)\to (S(Y),f(y_0))$ is small in $X_{k+1}$ due to attached space $SO_Y(\NN)_{f}$. \hfill $\blacksquare$

\begin{prop}
The natural inclusion $f\colon (Y,y_0)\to (Y_{s_0},y_0)\subset (X_0,y_0) \subset (S(Y),y_0)$ is homotopically non-trivial in $S(Y)$ if $Y$ is compact Hausdorff and homotopically non-trivial.
\end{prop}

\proof Suppose there exists a homotopy $H$ between $f$ and a constant map. By Lemma \ref{FactorHomotopy}   the homotopy $H$ factors over $\bigcup_{s\in T}W_s \cup \bigcup_{s\in S-T}Y_s$ for some finite subset $T\subset S$. Observe that  there exists a retraction of $\bigcup_{s\in T}W_s \cup \bigcup_{s\in S-T}Y_s$  to $Y_{s_0}$: retract $W_s$ to $Y_{s_0}$ for $s\in T$ and contract the rest to $y_0$. The composition of $H$ with such retraction implies that $Y$ is homotopically trivial, a contradiction. \hfill $\blacksquare$

\begin{corr} Space $S(Y)$ is a small $Y-$space if $Y$ is compact Hausdorff and homotopically non-trivial.
\end{corr}

Space $S(Y)$ is a universal example of a small $Y$-space in the following way.

\begin{prop}
Suppose $f\colon (Y_{s_0},y_0)\to (X,x_0)$ is a map to a small $Y-$ space $(X,x_0)$ where $Y_{s_0}\subset X_0 \subset FS(Y)$. Then $f$ extends over $(FS(Y),y_0).$
\end{prop}

\section{Homotopical closeness}

Homotopical closeness is a concept which generalizes homotopical smallness. Its development is motivated by considering how close the two loops are in some space.
Roughly speaking, loop $\a$ is close to loop $\b\not\simeq \a$ if for each $\eps>0$  there exists a homotopic representative $\a_\eps$ of $\a$ so that $d(\a_\eps(t),\b(t))<\eps, \forall t.$ In other words, there is no pair of close loops if the following condition holds: whenever there are homotopic loops $\a_\eps$ so that $\a_\eps(t)\stackrel{\eps\to 0}{\longrightarrow} \a(t)$ then $\a_\eps \simeq \a$. This condition is related to the property of being $\pi_1-$shape injective (or just shape injective) due to \cite{CC1} and to the property of being homotopically path Hausdorff.

\subsection{On homotopical smallness and closeness}

The aim of this section is to discuss some issues concerning the relationship between homotopical smallness and closeness. The following is the definition of closeness we employ for the future use.

\begin{dfn}\label{closeness}
Let $A\subset Y$ be a closed subspace of $Y$ and let $(X,d)$ be a metric space. Map $f\colon Y \to X$ is (homotopically) \textbf{close} to the map $g\colon Y \to X$ relatively to $A$ (denoted by $\rel A$) if the following conditions hold:
 \begin{itemize}
  \item [(a)] $f$ is not homotopic to $g, \rel A$ (i.e. there exists no homotopy between $f$ and $g$ that fixes all points of $A$);
  \item  [(b)] for each $\eps>0$ there exists a homotopy $H_\eps \colon Y \times [0,1]\to X$ so that
\begin{enumerate}
  \item $H_\eps|_{Y \times \{0\}}=f$;
  \item $H_\eps(a,t)=g(a), \quad\forall a\in A, \quad \forall t\in [0,1];$
  \item $d\big(H_\eps(y,1),g(y))< \eps, \quad\forall y\in Y.$
\end{enumerate}
\end{itemize}
\end{dfn}

The first observation is that (homotopical) closeness is only considered in metric spaces. The reason is that, roughly speaking,  we want to obtain homotopically equivalent maps  $f_n\colon Y \to X$ that point-wise uniformly converge to a map $f\colon Y \to X$ which is not homotopically equivalent to any $f_n$. The structure of a metric space was not required in the case of smallness (i.e. closeness to a constant map) as we were only considering convergence  towards one point. On the other hand the definition of closeness contains convergence of sequences with potentially different limit points. This generalization also allows us to consider closeness relatively to a subset $A$. Smallness could only be considered relatively to $\emptyset$ or a basepoint (yielding pointed and unpointed smallness).

\begin{figure}
\begin{center}
\scalebox{0.9}{\includegraphics{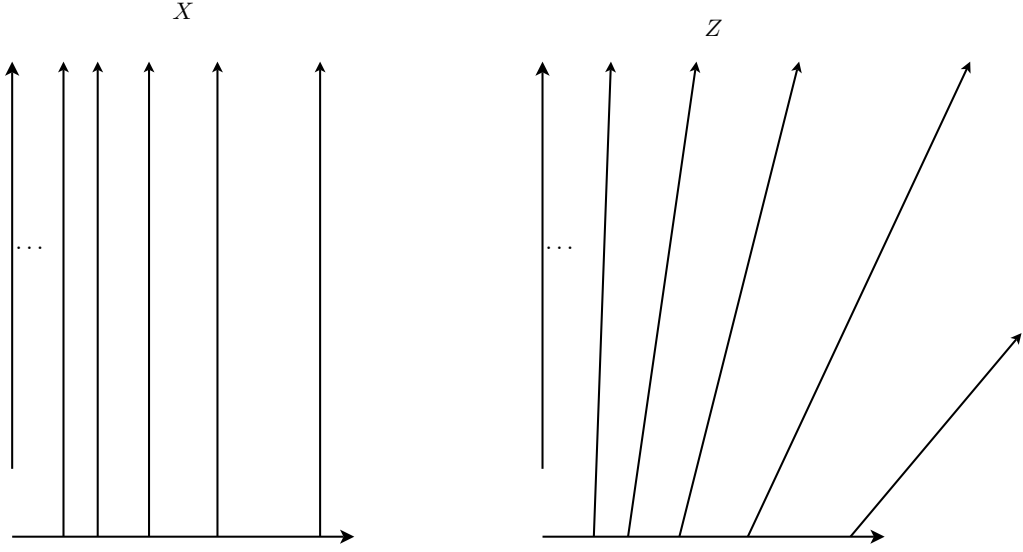}}
\end{center}
\caption{Spaces used to disprove the topological invariance of closeness.}
\label{closeness1}
\end{figure}

Another issue is the invariance of closeness and smallness. Every small map $Y\to X$ is topologically invariant (i.e. smallness is preserved by homeomorphisms of $X$). On the other hand the closeness is not preserved by homeomorphisms as suggested by the following example. Consider planar spaces $X$ and $Z$ (see Figure \ref{closeness1}) defined by the following rule.
$$
X:= \{x>0,y=0\} \cup \{x=0,y>1\} \cup \bigcup_{n\in \ZZ_+} \{x=\frac{1}{n},y>0\}
$$
$$
Z:= \{x>0,y=0\} \cup \{x=0,y>1\} \cup \bigcup_{n\in \ZZ_+} \{y=nx- 1; x> \frac{1}{n}\}
$$

Observe that there exists a homeomorphism $h\colon X\to Z$ which fixes the subset $\{x>0,y=0\} \cup \{x=0,y>1\}$ and linearly maps $\{x=\frac{1}{n},y>0\}$ to $\{y=nx- 1; x> \frac{1}{n}\}$. Consider the map $f\colon \RR_+\to X$ defined by $f(t):=(1,1+t)$. Note that $f$ is close (but not homotopic) to the map $g$ defined by $g(t):=(0,1+t)$. On the other hand $hf$ is not close to $hg$.

However, closeness is Lipschitz invariant and closeness of maps in a compact space is topologically invariant as proved by the following statements.

\begin{prop}
Let $f\colon Y \to X$ be close to $g\colon Y\to X$ (in a metric space $X$) and suppose a map $h\colon X \to Z$ is uniformly continuous. Then $hf$ is either close or homotopic to $hg$.
\end{prop}
\proof
For every $\eps>0$ let $\d_\eps$ denote a positive number so that if $d_X(x,y)<\d_\eps$ then $d_Z(f(x),f(y))<\eps.$ Consider homotopies $H_\eps$ according to Definition \ref{closeness}. Given a homeomorphism $h$ the homotopies $\widetilde H_\eps := hH_{\d_\eps}$ satisfy condition $(ii)$ of Definition \ref{closeness}.
\hfill$\blacksquare$

\begin{cor}\label{invariantCloseness}
Let $f\colon Y \to X$ be close to $g\colon Y\to X$ (in a metric space $X$) and let  $h\colon X \to Z$ be a map. If $X$ is compact Hausdorff or $h$ is Lipschitz then $hf$ is either close or homotopic to $hg$.
\end{cor}

Another observation is related to the nature of closeness as a relation, i.e. the absence of symmetry. Note that  $f$ being close to $g$ does not imply $g$ being close to $f$. The reason is that the definition of closeness ($f$ being close to $g$) requires homotopic representative of $f$ converging to the map $g$ and not only to  homotopic representative of $g$. The relaxation of condition of closeness as suggested by the last sentence would yield a symmetric (and transitive) relation of closeness. However, such relaxation would change the nature of closeness (and smallness) drastically. In particular,  every two maps to a graph of function $f(x)=x^{-2}$ would be either homotopic or close. Such relaxation would redefine small loops (within the unpointed category) in the following way: a loop is small iff it has a homotopic representative of diameter at most $\eps$ for every $\eps > 0$. This would mean that the punctured open disc is a small loop space but an open annulus is not hence smallness would not be a topological invariant. Also, the absence of close paths relatively to the endpoints would not coincide with the concept of the property of being homotopically path Hausdorff as proved by Proposition \ref{ClosePathsAndHomotopPathHaus}. For these reasons the definition of closeness is not symmetric.

The notion of small loops (in the pointed category) is closely related to the property of being homotopically Hausdorff (i.e. it is equivalent to the absence of small loops). In a similar fashion the closeness of paths relatively to the endpoints is related to the property of being homotopically path-Hausdorff in a locally path connected space.

\begin{prop}\label{ClosePathsAndHomotopPathHaus}
A locally path connected metric space $X$ has the property of being homotopically path Hausdorff iff there are no close paths $[0,1]\to X$ relatively to the endpoints of the interval.
\end{prop}

\proof
Suppose $X$ is not homotopically path Hausdorff. According to Definition \ref{homotopHauss} there exist paths $w,v\colon [0,1]\to X$, $v \not\simeq w \rel \{0,1\}$ with the following property: for any chosen $n\in \ZZ_+$ and a cover of $w([0,1])$ by open sets of diameter at most $\frac{1}{n}$ the conditions of Definition \ref{homotopHauss} are not satisfied due to some path $v_n$ homotopic to $v$ relatively to $\{0,1\}$. In particular,  $d(w(t),v_n(t))<\frac{1}{n}$, $\forall t\in [0,1]$ as $w(t)$ and $v_n(t)$ are contained in a set of diameter at most $\frac{1}{n}$. This implies that $v$ is close to $w$ relatively to $\{0,1\}$.

To prove the other direction we use local path connectedness of $X$. Suppose the path $v\colon [0,1]\to X$ is close to the path $w$ relatively to $\{0,1\}$. Choose any cover of $w([0,1])$ by finitely many path connected open sets $U(P_1),\ldots, U(P_k)$ and any partition $0=t_0 < t_1 < t_2 < t_3 <\ldots< t_k=1$ so that:
\begin{enumerate}
  \item $U(P_j)$ covers $w([t_{j-1},t_j])$;
  \item $P_j\in w([t_{j-1},t_j])$.
\end{enumerate}
There exists $\eps>0$ so that the $\eps$-neighborhood of $w([t_{j-1},t_j])$ is contained in $U(P_j), \forall j$. The closeness of $v$ to $w$ relatively to $\{0,1\}$ allows us to choose a path $v'\simeq v \rel \{0,1\}$ so that $d_X(v'(t),w(t))< \eps$ hence $U(P_j)$ covers $v'([t_{j-1},t_j])$. Since the sets $U(P_j)$ are path connected we can (for each $j$) redefine $v'|_{[t_{j-1},t_j]}$ so that we do not change the homotopy type relatively to $\{t_{j-1},t_j\}$, $P_j\in v'([t_{j-1},t_j])$ and $ U(P_j)\supset v'([t_{j-1},t_j])$. Such path $v'$ contradicts  Definition \ref{homotopHauss} for the loop $\a=w*v^-$ hence $X$ is not homotopically path Hausdorff. \hfill$\blacksquare$
\bigskip

If the space $X$ is not locally path connected then close paths need not contradict the property of being homotopically path Hausdorff. This fact is connected to the following observation. If a path $f\colon [0,1]\to X$ is close to the path $g$ (relatively to $\{0,1\}$) then:
\begin{itemize}
  \item [(i)] $f$ is close to $g$ in the Peanification $PX$  if $X$ is locally path connected.
  \item [(ii)] $f$ may not be close to $g$ in the Peanification $PX$ (for some metric on $PX$) if $X$  is not locally path connected.
\end{itemize}

Statement $(i)$ is obvious as $X=PX$ in the case of a locally path connected space. Note that statement $(ii)$ requires a structure of metric space on $PX$ in order to consider closeness. To prove statement $(ii)$ we construct the space $C(S^1,\{0\})$ \index{C$(S^1,\{0\})$ space} which is a modification of $HA$. Recall that $HA$ is constructed with the aim to create a small loop. In order to do this we attach big homotopies along the loops converging to a point. The construction of $C(S^1,\{0\})$ follows the same philosophy for closeness. We attach big homotopies along loops converging to another loop (rather than a point). Recall that $0=(0,0)\in \RR^2$ and $S^1(S,r)$ denotes a circle in $\RR^2$ with center $S$ and radius $r$. Define
$$
S^1_n:=S^1((1+\frac{1}{n},1),1+\frac{1}{n}) \quad \emph{ for } n\in \ZZ_+,\quad
S^1_\infty :=S^1((1,0),1).
$$

Naturally embed $\cup_{n\in \ZZ_+}S^1_n\cup S^1_\infty$ in $\RR^3$ and attach spaces $A_n = (S^1\times [0,1])_n$ for all $n\in \ZZ_+$ so that:
\begin{enumerate}
  \item we identify $(S^1\times \{0\})_n$ with $S^1_n$;
  \item we identify $(S^1\times \{1\})_n$ with $S^1_{n+1}$;
  \item we identify $(\{0\}\times [0,1])_n$ with $0$;
  \item the rest of $A_n$ is stretched up so that it reaches height ($z-$coordinate) $1$ (i.e. $A_n$ has a point which is of distance at least $1$ from $S^1_{n}$ and $S^1_{n+1}$).
\end{enumerate}

\begin{figure}
\begin{center}
\includegraphics{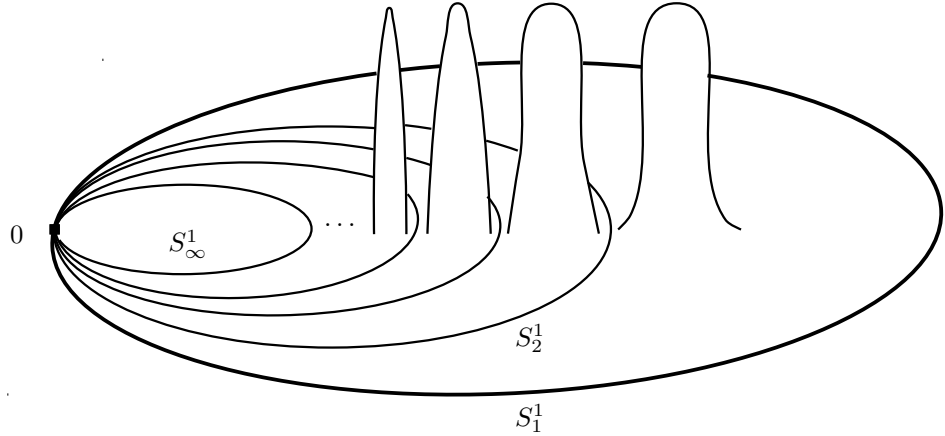}
\end{center}
\caption{The space $C(S^1,\{0\})$.}
\label{CloseHA}
\end{figure}

In other words, we attach big homotopies between loops $S^1_n$ as suggested by Figure \ref{CloseHA}. Since closeness is only defined in a metric space we can not attach all $A_n$ by quotient maps (as in $SHA$) but rather within a metric space $\RR^3$.
\bigskip

\remark The conditions above about the nature of the attached homotopies $A_n$ do not uniquely define the space $C(S^1,\{0\})$. The reason is, roughly speaking, that the homotopies $A_n$ approach the loop $S^1_\infty$ rather than just a point. For example, the homotopies $A_n$ may be chosen so that for any given $x\in S^1_\infty-\{0\}$ the space either is or is not locally path connected at $x$. Figure \ref{CloseHA} suggests that $C(S^1,\{0\})$ is locally path connected everywhere and that there are small loops at $(1,0,0)\in S^1_\infty$ as the humps (i.e. subspaces of $A_n$ with $z-$coordinate at least $1$) of $A_n$ converge to that point. In order to comply with later definitions we demand the humps  of $C(S^1,\{0\})$ to converge to entire $S^1_\infty$ so that $C(S^1,\{0\})$ is not locally path connected at any point of $S^1_\infty -\{0\}$. For an alternative description of $C(S^1,\{0\})$ see Definition \ref{CodYAspace}. \bigskip

Note that the loops $S^1_n$ are homotopic to each other relatively to $0$ via homotopies $A_n$ but these homotopies can not be combined to obtain a homotopy to the limit loop $S^1_\infty$. Similarly as in the case of $HA$ we can prove that the map $(S^1,0)\to (S^1_1,0)\subset C(S^1,\{0\})$ is close to the map $(S^1,0)\to (S^1_\infty,0)\subset C(S^1,\{0\})$ relatively to $0$. In both cases we consider a map of the form $e^{i\f} \mapsto (A+Be^{i\f},0).$ However, the Peanification of $C(S^1,\{0\})$ is a wedge of $S^1$ and $\big( S^1\times [0,1) \cup \{0\}\times \{1\}\big)$ hence it contains no close loops. The space $C(S^1,\{0\})$ is an example of a homotopically path Hausdorff space with close loops.

The last observation is related to the Spanier group of a space. Groups $\pi^s$ and $\pi^{sg}$ are generated by small loops and defined in \cite{Vi1}. Close loops have no influence on these groups but  may interfere with the Spanier group. The following is the definition of a Spanier group for locally path connected spaces as presented in  \cite{4aut}.

\begin{dfn}\label{SpanierGrp}
Let $(X,x_0)$ be a locally path connected space and let ${\mathcal U} = \{U_i\}_{ i \in I }$ be a cover of $X$ by open neighborhoods. Define $\pi_1({\mathcal U}, x_{0})$ as the subgroup of $\pi_1(X,x_0)$ consisting of the homotopy classes of loops that can be represented by a product (concatenation) of the following type:
$$
\prod^n_{j=1} u_j*v_j*u_j^{-},
$$
where  $u_j$ are paths that run from $x_0$ to a point in some $U_i$ and each $v_j$ is a closed path inside the corresponding $U_i$ based at the endpoint of $u_j$. We call $\pi_1({\mathcal U}, x_{0})$ the \textbf{Spanier group of $(X,x_0)$ with respect to ${\mathcal U}$}.

Let ${\mathcal U}$  and ${\mathcal V}$ be an open covers of $X$ and let ${\mathcal U}$ be a
refinement of ${\mathcal V}$. Then $\pi_1({\mathcal U},x_0) \subset \pi_1({\mathcal V},x_0)$.
This inclusion relation induces an inverse limit defined via the directed system of all covers with respect to refinement. We will call such limit the  \textbf{Spanier group} of the space $X$ and denote it by $\pi^{sp}_1(X,x_0)$.
\end{dfn}

\begin{prop} Let $(X,x_0)$ be a locally path connected space.
\begin{enumerate}
  \item $\pi^{sg}_1(X,x_0) \subset \pi^{sp}_1(X,x_0)$.
  \item If $f\colon ([0,1],0)\to (X,x_0)$ is close to $g$ relatively to $\{0,1\}$ in a metric space $(X,d)$ then $[f*g^-]\in \pi^{sp}_1(X,x_0)$.
\end{enumerate}
\end{prop}

\proof
Claim (i) is true as every element of $\pi^{sg}_1(X,x_0)$ is contained in each $\pi_1({\mathcal U}, x_{0})$ by the definition.

To prove claim (ii) we partially imitate the proof od \ref{ClosePathsAndHomotopPathHaus}. Fix a cover ${\mathcal U}$ of $X$ and choose a finite subfamily $U_1,\ldots,U_k \subset {\mathcal U}$ covering $g([0,1])$ so that for some partition $0=t_0<t_1<\ldots<t_k=1$ the set $U_j$ contains $g([t_{j-1},t_j]), \forall j.$ There exists $\eps > 0$ so that for every $j$:
\begin{itemize}
  \item $\eps$-neighborhood of $g([t_{j-1},t_j])$ is contained in $U_j$;
  \item every point of the $\eps-$neighborhood of $g(t_j)$ is connected to $g(t_j)$ by a path in $U_j\cap U_{j+1}$.
\end{itemize}

\begin{figure}
\begin{center}
\includegraphics{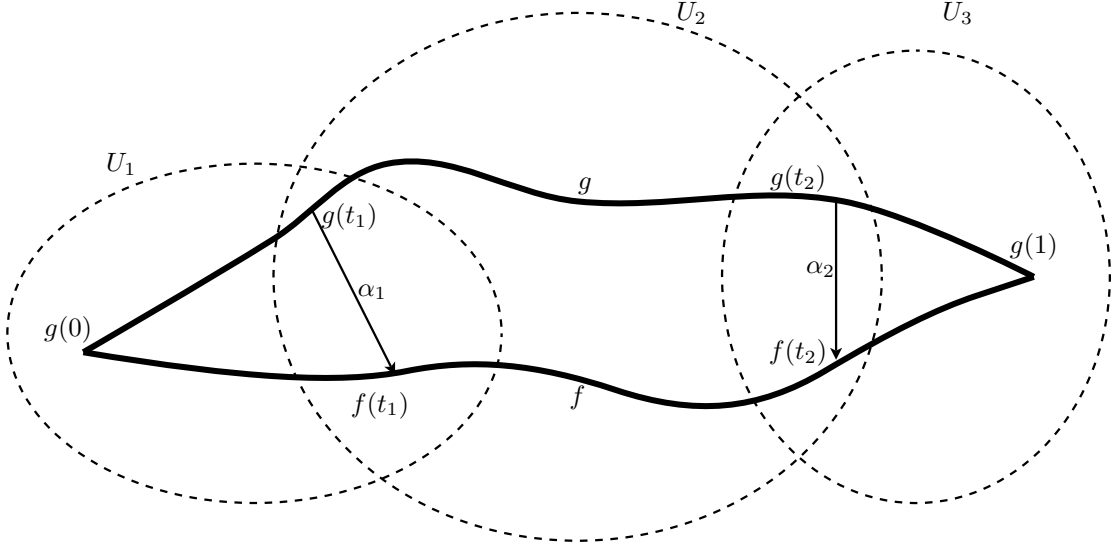}
\end{center}
\caption{Notation concerning close paths $f$ and $g$.}
\label{spanier}
\end{figure}

We can assume $d_X(f(t),g(t))<\eps, \quad \forall t$. For each $j$ let $\a_j$ denote an oriented path in $U_j\cap U_{j+1}$ between $g(t_j)$ and $f(t_j)$ as denoted by Figure \ref{spanier}. We can assume $\a_0$ and $\a_k$ to be constant paths.  Observe that the oriented loop $Q_j$ defined as a concatenation

$$
\a_{j-1} * f|_{[t_{j-1},t_{j}]} * \a_j^- * \big(g|_{[t_{j-1},t_{j}]}\big)^-
$$

is based at $g(t_j)$ and contained in $U_j$.
The class $[f*g^-]$ is contained in $\pi_1({\mathcal U},x_0)$ because it can be expressed as
$$
\prod^k_{j=1} g|_{[0,t_{j-1}]} * Q_j  * \big(g|_{[0,t_{j-1}]}\big)^-.
$$
Hence $[f*g^-]\in \pi^{sp}_1(X,x_0)$.
\hfill$\blacksquare$

\subsection{Standard constructions}

In this subsection we present some aspects of closeness which are motivated by similar results on smallness. Since the closeness is only defined in metric spaces the construction of an m-stratified space and some other features of smallness are not applicable. The absence of these obstruct the generalization of some constructions including the small loop space. However it is possible to construct the space $C(Y,A)$ with maps $f,g\colon Y\to C(Y,A)$ for which the map $f$ is close to $g$ relatively to $A$.

Given a metric space $X$ its metric will be denoted by $d$ or $d_X$. The metric on a product of metric spaces is defined by $d_{X\times Y}\big((x_1,y_1),(x_2,y_2)\big):=d_X(x_1,x_2)+d_Y(y_1,y_2).$

\subsubsection{Free closeness}

We first consider free closeness, i.e. closeness relatively to $\emptyset$. The following is a generalization of the $SO$ spaces.

\begin{dfn} \label{CodYspace}
Let $Y$ be a metric space. Metric space $C(Y)$ is the subspace of $Y\times [0,1]\times [-1,1]$ defined as
$$
\Big\{(y,t,\sin \frac{\pi}{t}); (y,t)\in Y \times (0,1]\Big\}\quad \bigcup \quad Y \times \{0\}\times \{0\}.
$$
\end{dfn}

\begin{figure}
\begin{center}
\includegraphics{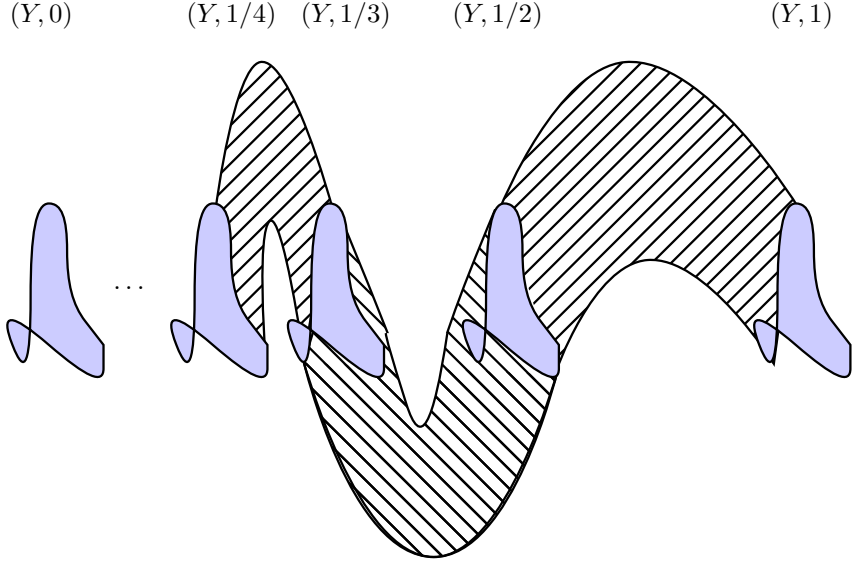}
\end{center}
\caption{Space $C(Y).$}
\label{C(Y)}
\end{figure}

Figure $\ref{C(Y)}$ schematically represents space $C(Y)$: copies of space $Y$ connected by big (dashed) homotopies converge to $(Y,0,0)$. In the case of $Y$ being a single point we obtain  $C(Y)=\{(0,0)\}\cup \{(x,\sin\frac{1}{x}); x\in (0,1]\}$.
For any point $(y,t,s)\in C(Y)$ we refer to $y,t,s$ as the first, the second and third coordinate respectively. The role of these coordinates is the following:
\begin{itemize}
  \item the first coordinate allows space $Y$ to be embedded;
  \item the second coordinate represents homotopies between converging embeddings;
  \item the third coordinate makes homotopies via the second coordinate big so that they cannot extend over $(Y,0,0)$.
\end{itemize}

Note that $(Y,0,0)$ is not path connected to $(Y,1,0)$ which yields the following result.

\begin{prop}\label{BB}
The map $f\colon Y\to C(Y)$ defined by $y\mapsto (y,1,0)$ is close to the map $g\colon Y\to C(Y)$ defined by $y\mapsto (y,0,0)$
\end{prop}

\begin{prop}\label{BBB}
Suppose  $Y$ is a compact metric space and the map $f\colon Y \to X$ is not homotopic to $g\colon Y \to X$. The map $f$ is close to $g$ iff there exists a map $F\colon C(Y)\to X$ so that $F|_{(Y,1,0) }=f$ and $F|_{(Y,0,0)}=g$.
\end{prop}

\proof
The existence of an extension $F$ implies that $f$ is close to $g$ by Corollary \ref{invariantCloseness}.

To prove the other direction assume that $f$ is close to $g$. Hence for all $n\in \ZZ_+$ there exist maps
$$
H_n\colon Y \times \Big[\frac{1}{n+1},\frac{1}{n}\Big] \to X; \quad H_n|_{Y\times\{n+1\}}=H_{n+1}|_{Y\times\{n+1\}};
$$
$$
d\Big(H_n (y,\frac{1}{n}),g(y)\Big) < \frac{1}{n}; \quad H_1|_{Y \times \{1\}}=f.
$$
We have to adjust the maps $H_n$ in order to construct a continuous map on $C(Y)$. The idea is to adjust the maps $H_n$ so that they only depend on $Y-$coordinate in   appropriate neighborhoods of $(Y,n+1)$ and $(Y,n)$.

Given any $n\in \ZZ_+$ and any $\d < (\frac{1}{n}-\frac{1}{n+1})/2$ we can assume $H_n(y,t)=H_n(y,\frac{1}{n+1})$ if $|t-\frac{1}{n+1}|< \d$ and $H_n(y,t)=H_n(y,\frac{1}{n})$ if $|t-\frac{1}{n}|< \d$. A required modification can be obtained as follows. Extend $H_n$ to a map $Y \times [a,b] \to X$ so that $$
H_n(y,t):=H_n(y, \frac{1}{n+1}) \quad \emph{ if } \quad  t < \frac{1}{n+1};
$$
$$
H_n(y,t):=H_n(y, \frac{1}{n}) \quad \emph{ if } \quad  t > \frac{1}{n}.
$$
The linear contraction $c\colon [a,b]\to [\frac{1}{n+1},\frac{1}{n}]$ for appropriate $a$ and $b$ induces a map
$$
Y \times \Big[\frac{1}{n+1},\frac{1}{n}\Big] \stackrel{1\times c\-}{\longrightarrow} Y \times [a,b] \stackrel{H_n}{\longrightarrow} X
$$
which satisfies the required condition.

Define the map $F\colon C(Y)\to X$ by the rule $(y,t,s)\mapsto H_n(y,t)$ if $t\in [\frac{1}{n+1},\frac{1}{n}]$ and  $F(y,0,0)=f(y)$.  We claim that for a suitable choice of maps $H_n$ the map $F$ is continuous. Note that $C(Y)\cap \big(Y\times (0,1]\times [-1/2,1/2]\big)$ is a disjoint union of closed neighborhoods of $(Y,1/n,0)$  which are all homeomorphic to $Y \times [-1,1]$. We can assume (by applying the modification above) that the maps $H_n$ are appropriately modified to ensure that $F(y,s,t)$ only depends on $y$ on each of these  sets. In particular, any sequence $a_n=F(y,t_i,s(t_i))$ with fixed $y$ and $t_i$ converging to $1/n$ is eventually  constant. The stabilization occurs (if not before) for $i$ with the property that for all successive indexes $j>i$ we have  $|s(t_j)|< 1/2$ and $|t_j-1/n|<1/(4n^2)$.

Note that $F$ is continuous at every point $(y,t,s)$ with $t>0$ as the maps $H_n$ are continuous and agree on the intersection of their domains. To prove that $F$ is continuous consider a convergent sequence $(y_i,t_i,s_i)\to (y_0,0,0)$ in $C(Y)$. Given any $\eps>0$ choose $i_0$ so that for every $i>i_0$:
\begin{itemize}
  \item $d\big(f(y_i),f(y_0)\big)< \eps/2$;
  \item $|s_i|< 1/2$;
  \item $t_i < 1/n_\eps<\eps/2$ for some $n_\eps\in \ZZ_+$ (i.e. $d\big( f(y,t_i,s_i),f(y,0,0) \big)< \frac{1}{n_\eps}, \forall y\in Y$).
\end{itemize}

Then
$$
d\big(F(y_i,t_i,s_i),F(y_0,0,0)\big)<
$$
$$
<d\big(F(y_i,t_i,s_i),F(y_i,0,0)\big)+d\big(F(y_i,0,0),F(y_0,0,0)\big)<1/n_\eps + \eps/2 < \eps,
$$
since $F(y_i,t_i,s_i)\in \{F(y_i,1/n,0)\}_{n \geq n_\eps}$ hence $F$ is continuous.\hfill$\blacksquare$

\subsubsection{Relative closeness}

\begin{dfn}\label{CodYAspace}
Let $Y$ be a metric space and let $A\subset Y$ be a closed subspace. Choose a map $\f\colon Y \to [0,1]$ so that $A=\f \i (\{0\})$. The metric space $C(Y,A)$ is a subspace of $Y\times [0,1]\times [-1,1]$ defined as
$$
\Big\{(y,\f(y)t,\f(y)\sin \frac{\pi}{t}); (y,t)\in (Y\backslash A) \times (0,1]\Big\}\quad \bigcup \quad Y \times \{0\}\times \{0\}.
$$
\end{dfn}

The space $C(Y,A)$ depends on a choice of  map $\f$. Nevertheless we omit $\f$ form the notation of $C(Y,A)$ as the properties of our interest do not depend on the choice of a map $\f$. Note that $C(Y,A)$ is not locally path connected at any point of $(Y\backslash A,0,0)$. In the case of $(Y,A)=([0,1],\{1\})$ the space $C(Y,A)$ is a cone over the space
$$
C(\{0\})=\{(0,0)\}\cup \{(x,\sin\frac{1}{x}); x\in (0,1]\}.
$$
The following proposition provides some examples of a relatively close maps.

\begin{prop}\label{B}
Suppose $A\subset Y$ is a closed subspace of a metric space $Y$. The inclusion $i_1\colon Y \hookrightarrow C(Y,A)$ defined by $y \mapsto (y,\f(y),\f(y))$ is homotopic or close to the inclusion $i_2\colon Y \hookrightarrow C(Y,A)$ $\rel A$, where  $i_2\colon y\mapsto (y,0,0)$.
\end{prop}

\proof
The homotopies $H\colon Y \times [1/n,1]\to C(Y,A)$ defined by the rule
$$
(y,t)\mapsto (y,\f(y)t,\f(y)\sin \frac{\pi}{t})
$$
homotope $i_1$ arbitrarily close to $i_2$. \hfill$\blacksquare$
\bigskip

In the case of free closeness the inclusions  of Proposition \ref{BB} were not homotopic due to an argument on path connectedness. Such argument can not be employed  for the inclusions $i_1$ and $i_2$ of Proposition \ref{B} as $C(Y,A)$ is path connected if $Y$ is path connected and $A \neq \emptyset$. In order to find a condition for the inclusions $i_1$ and $i_2$ not to be homotopic relatively to $A$ consider the space $W:=Y_1 \cup_{a\sim (a,(0,1]);  a\in A} \big(Y_2 \times (0,1]\big)$ where $Y_1$  and $Y_2$ are isomorphic to a locally path connected metric space $Y$. Using the argument of  Peanification for the points of $(Y-A,0,0)$ observe that the natural bijection $W \leftrightarrow C(Y,A)$ induces a natural bijection on maps and homotopies from $Y, \rel A$ if $Y$ is locally path connected. The space $W$ is homotopic to $Y_1\cup_{A} Y_2$ with the natural inclusions of $Y_1$ and $Y_2$ corresponding to  the inclusions $i_1$ and $i_2$ of Proposition \ref{B}. Hence the inclusions $i_1$ and $i_2$ of Proposition \ref{B} are homotopic $\rel A$ iff the inclusions of $Y_1$ and  $Y_2$ into $Y_1\cup_{A} Y_2$  are homotopic $\rel A$. The later of these conditions is equivalent (via contraction of $Y_2$ in $Y_1\cup_{A} Y_2$ ) to $Y/A$ being contractible which is the case if $A \hookrightarrow X$ is a cofibration and a homotopy equivalence.

\begin{cor}\label{BBBB}
Suppose the cofibration inclusion $A\hookrightarrow Y$ of a closed subspace into a locally path connected metric space $Y$ is not a homotopy equivalence. Then the inclusion $i_1\colon Y \hookrightarrow C(Y,A)$ defined by $y \mapsto (y,\f(y),\f(y))$ is close to the inclusion $i_2\colon Y \hookrightarrow C(Y,A)$ $\rel A$, where  $i_2\colon y\mapsto (y,0,0)$.
\end{cor}

Combining  Corollary \ref{BBBB} and Proposition \ref{BBB} we obtain the extensional classification of certain types of relative close maps.

\begin{prop}
Consider the following situation:
\begin{itemize}
  \item $Y$ is a compact, locally path connected metric space;
  \item $A \subset Y$  is a closed subspace;
  \item the natural inclusion $A \hookrightarrow Y$ is a cofibration which is not a homotopy equivalence;
  \item the map $f\colon Y \to X$ is not homotopic to the map $g\colon Y \to X, \rel A$;
  \item the inclusions $i_1,i_2\colon Y \to C(Y,A)$ are defined by $i_1(y)=(y,\f(y),\f(y))$ and $i_2(y)=(y,0,0).$
\end{itemize}
Then the map $f$ is close to $g, \rel A$ iff there exists a map $F\colon C(Y,A)\to X$ so that $F i_1=f$ and $F i_2=g$.
\end{prop}

\proof
The existence of an extension $F$ implies that $f$ is close to $g$ by Corollary \ref{invariantCloseness}. To prove the other direction assume that $f$ is close to $g, \rel A$. By Proposition \ref{BBB} there exists an appropriate extension $F$ over $C(Y)$. Note that for every $a\in A$ the closed subset $W_a:=C(Y)\cap \big(\{a\}\times[0,1]\times [-1,1]\big)$ is mapped by $F$ to $a$. The map $F$ induces an extension over the quotient space $C(Y)/_{W_a; a\in A}=C(Y,A)$.\hfill$\blacksquare$

\subsection{Closeness and compactness}

Closeness of maps due to Definition \ref{closeness} is considered within metric spaces in order to enforce a uniform continuouity of the approaching maps. On the other hand the idea of close maps appears in (iii) of Definition \ref{homotopHauss} where closeness of paths is considered in a non-metric space. The aim of this section is to introduce the notion of closeness for maps with compact domain and possibly non-metric range.

\begin{dfn}\label{CloseINCompact}
Suppose $f,g\colon K \to Y$ are maps defined on a compact Hausdorff space $K$ so that $f\not\simeq g,\rel A$ for some closed subspace $A\subset K$. The map $f$ is close to $g$ if for every finite open cover $U_1, \ldots, U_k$ of $g(K)$ there exist:
 \begin{itemize}
   \item  a collection $B_1, \ldots, B_k$ of closed subsets of $K$ so that $K=\cup_i B_i$ and $g(B_i) \subset U_i$;
   \item the map $f' \simeq f, \rel A$ so that $f'(B_i)\subset U_i, \forall i.$
 \end{itemize}
\end{dfn}

Let us prove that both definitions of closeness agree if considered in a metric space $Y$. Suppose the map $f\colon K \to Y$ is close $\rel A$ to the map $g$ in terms of Definition \ref{CloseINCompact} where $A\subset K$ is a close subset of a compact Hausdorff space $K$ and $Y$ is a metric space. Given any $\eps>0$ we can cover $g(K)$ by a collection of open sets of diameter at most $\eps$. The map $f'$ referred to such cover by Definition \ref{CloseINCompact} is homotopic to $f, \rel A$ and satisfies the condition $d(f'(x),g(x))<\eps, \forall x\in K$, hence $f$ is close to $g$ in terms of Definition \ref{closeness}.

To prove the opposite implication assume that the map $f\colon K \to Y$ is close $\rel A$ to the map $g$ in terms of Definition \ref{closeness}. Given any finite open cover $U_1, \ldots, U_k$ of $g(K)$ choose a collection $B_1, \ldots, B_k$ of closed subsets of $K$ so that $K=\cup_i B_i$ and $B_i \subset U_i$. There exists an $\eps>0$ so that for all $i$ the $\eps-$neighborhood of $B_i$ is also in $U_i$. A map $f'\simeq f, \rel A$ with the property of $d(f'(x),g(x))<\eps, \forall x\in K$ satisfies the conditions of Definition \ref{CloseINCompact} hence the definitions are equivalent.

\section{Applications}

Other that the classification of Theorem \ref{classifOFhomotopt2}, homotopical smallness and closeness can efficiently be used in construction of certain spaces. The simplest case is the topologist's sine curve, which is equivalent to $C(\{0\}\})$. It can also be considered as a one-dimensional harmonic archipelago $HA^1$. The harmonic archipelago of dimension $n$ [denoted by $HA^n$] is a subset of $\RR^{n+1}$ is constructed from a wedge of spheres
$\{S^n_i\}_{i\in \ZZ_+}$ radii $1/i$ by attaching big homotopies (i.e. of diameter at least 1) between each pair of consecutive spheres $(S^n_i, S^n_{ i+1} )$.

\begin{figure}
\includegraphics[scale=0.60]{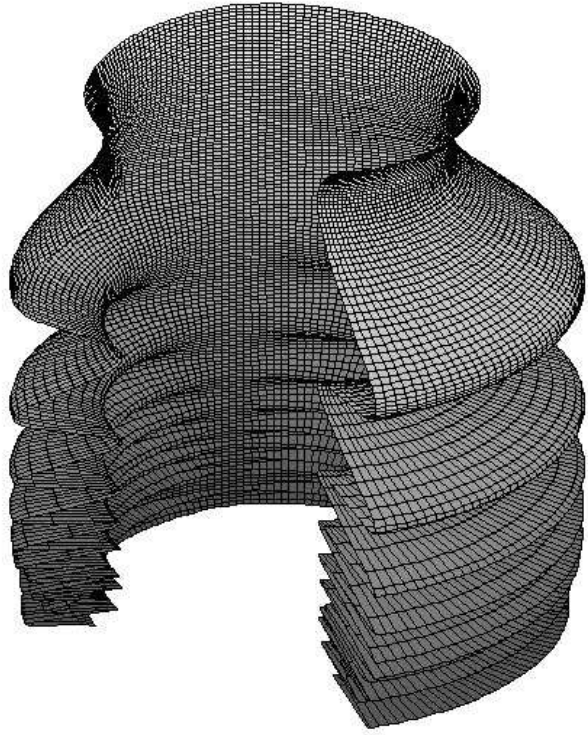}
\caption{The surface portion of the harmonic vase.}
\label{HVase}
\end{figure}

\textbf{The harmonic vase}  was defined in \cite{Vi2} as a subset of $\RR^3$. It has an essential role in the proof of Theorem \ref{mojizrek}. It consists of a disc
$$
B^2=\{(x,y,0)\in \RR^3 ; x^2+y^2\leq 9\},
$$
and a surface portion (see Figure \ref{HVase})
$$
r:= \frac{|\f|}{\pi}\sin \frac{\pi}{z}+2, \quad  z\in (0,1],\quad \f\in [-\pi,\pi],
$$
where $(r,\f)$ are polar coordinates in $\RR^2\times \{0\}\subset \RR^3$ and $z$ is the  coordinate of $\{0\}^2\times \RR$ so that $(r,\f,z)$ are cylindric coordinates in $\RR^3$. The motivation for $HV$ is a construction of a loop $f$ which is close to a homotopically trivial embedding $g$ of $S^1$ in a compact space. The quotient space $HV/\{\f=0\}$ is homeomorphic to a compact space $C(S^1,\{0\})\cup B^2$ where the disc $B^2$ is attached in appropriate way as described above. Furthermore, the union of the surface portion of $HV$ and $\di B^2$ is equivalent  to $C(S^1,\{0\})$.

\begin{figure}
\includegraphics[scale = 0.7]{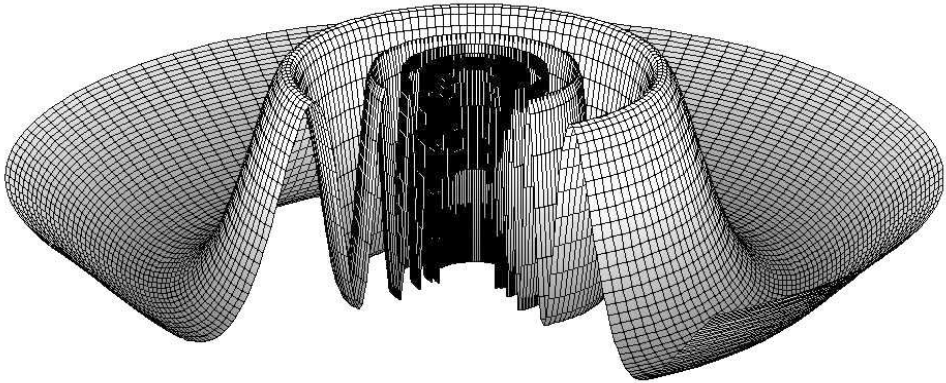}
\caption{The surface portion of space $A$.}
\label{spaceY}
\end{figure}

\textbf{Space A} as defined in \cite{Re} was developed as an example of a space which is homotopically Hausdorff but not strongly homotopically Hausdorff, i.e., it has no small loops but has free small loops. It is a subspace of $\RR^3$ consisting of three parts:
\begin{itemize}
  \item the surface portion which is obtained by rotating the topologist's sine curve
        $$
            \{(0,0,0)\}\cup \{(x,0,\sin\frac{1}{x}); x\in (0,1]\}
        $$
        around the $z-$axis, as suggested by Figure \ref{spaceY};
  \item the central limit arc $\{0\}\times \{0\}\times [-1,1]$;
  \item connecting arcs, i.e. a system of countably many closed radial arcs emerging from the central limit arc so that $A$ is compact and locally path connected.
\end{itemize}

The surface portion of $A$ can be considered as an unpointed version of $HA$, the quotient $C(S^1)/_{(S^1,0,0)}$ or as a modified version (in the same way as $HA$ is a modified version of $SO_{S^1}(\NN)$)  of $FSO_{S^1}(\NN)$.

Space $A$ and its properties  were studied in \cite{4aut} under the name of space $Y'$. A similar space called $Y$ is defined and studied in the same paper. It consists of the same surface portion and the same central limit arc, but instead of connecting arcs there is a single simple arc connecting the central limit arc with the surface portion. The space $Y$ is not locally path connected but the connecting arc makes it path connected.

Space $Y$ distinguishes between two definitions of semi-local simple connectedness that appear in the literature. According to \cite{4aut}:
 \begin{itemize}
   \item space $X$ is based semi-locally simply connected iff every point $x\in X$ has a neighborhood $U\subset X$ so that $\pi_1(U,x)\to \pi_1(X,x)$ is trivial.
   \item space $X$ is unbased semi-locally simply connected iff every point $x\in X$ has a neighborhood $U\subset X$ so that every loop in $U$ is contractible.
 \end{itemize}
Both definitions of semi-local simple connectedness agree if the space is locally path connected. It turns out that the space $Y$ is based but not unbased semi-locally simply connected due to its topology at the central limit arc.

\begin{figure}
\includegraphics[scale=1.1]{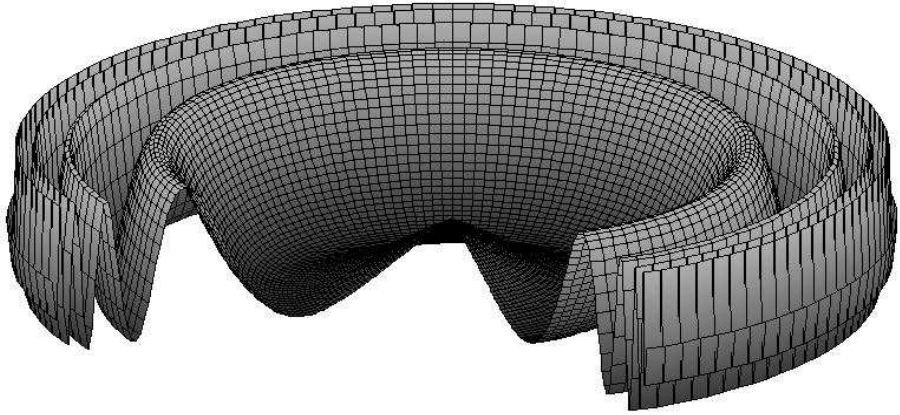}
\caption{The surface portion of space $B$.}
\label{spaceZ}
\end{figure}

\textbf{Space B} as defined in \cite{Re} was developed as an example of a space which is strongly homotopically Hausdorff but not shape injective. It is a subspace of $\RR^3$ consisting of three parts:
\begin{itemize}
  \item the surface portion which is obtained by rotating the topologist's sine curve
        $$
            \{(0,0,0)\}\cup \{(x,0,\sin\frac{1}{x}); x\in (0,1]\}
        $$
        around the axis $\{1\}\times \{0\}\times \RR$, as suggested by Figure \ref{spaceZ};
  \item the outer annulus obtained by rotating  $\{0\}\times \{0\}\times [-1,1]$ around the axis $\{1\}\times \{0\}\times \RR$;
  \item connecting arcs, i.e. a system of countably many closed radial arcs emerging from the outer annulus so that $B$ is compact and locally path connected.
\end{itemize}

The surface portion of $B$ is essentially the same as  $C(S^1)$. Space $B$ and its properties  were studied in \cite{4aut} under the name of space $Z'$. A similar space called $Z$ is defined and studied in the same paper.

\subsection{Realization theorems}

One of the basic problems in homotopy theory is the realization of various groups as a homotopy invariants of certain spaces. In particular, we are interested in the following question: given a group $G$ when can we realize it as a fundamental group of a path connected space $X$ which possesses the following properties:
\begin{enumerate}
  \item $X$ is compact;
  \item $X$ is metric;
  \item $X$ is locally path connected?
\end{enumerate}

It turns out that these three conditions are too restrictive for the realization of all countable groups.

\begin{thm}\cite{SH1}
Let $X$ be a compact metric space which is path connected and locally path connected. If the fundamental group of $X$ is not finitely generated then it has the power of the continuum.
\end{thm}

An improvement to this theorem has been made in \cite{CC1} and \cite{DV}.

\begin{thm}
Let $X$ be a compact metric space which is path connected and locally path connected. If the fundamental group of $X$ is not finitely presented then it has the power of the continuum.
\end{thm}

However, if we omit any of the three properties mentioned above we can realize all countable groups. It is well known that every group can be realized as a fundamental group of a path connected $CW$ complex of dimension two and every countable group can be realized as a fundamental group of a countable path connected $CW$ complex of dimension two. Every countable $CW$ complex  is  homotopy equivalent to a locally finite (hence metrizable) $CW$ complex of the same dimension which yields the realization in terms of metric, locally path connected spaces. Since the metric space of such realization is a two dimensional $CW$ complex it can be embedded in $\RR^5$.

\begin{thm}
Let $G$ be a countable group. Then $G$ can be realized as a fundamental group of a two dimensional metric space $X$ which is path connected and locally path connected.
\end{thm}

This result implies that, omitting the compactness from the list above, we can realize all countable groups as a fundamental groups of a space with prescribed properties. Similarly we can omit metrizability in order to obtain a realization in terms of compact locally path connected space.

\begin{thm}\cite{Paw}
Let $G$ be a countable group. Then $G$ can be realized as a fundamental group of a compact space $X$ which is path connected and locally path connected.
\end{thm}

The realization in terms of a compact metric space was proven in \cite{Vi2} using the techniques of homotopical closeness (i.e. the harmonic vase and its variation: the braided harmonic vase) and the universal Peano space. It turns out that given a locally path connected space $X$ in certain circumstances, one can construct a compact space $Y$ so that $PY\simeq X$, i.e. the spaces have the same fundamental group.

\begin{thm}\label{mojizrek}
Let $G$ be a countable group. Then $G$ can be realized as a fundamental group of a two dimensional compact metric space $X\subset \RR^4$ which is path connected.
\end{thm}

The approach of \cite{Vi2} in terms of homotopical smallness was generalized in \cite{TeVi} in order to obtain a wider class of realization theorems. This improvement includes the realization of appropriately prescribed groups as a  homotopy or homology  groups of a space. The realization results are implied by the following fact.

\begin{prop}\cite{TeVi}\label{TeVi}
For every countable $CW$ complex $K$ there is a compact metric space $X$ such that $PX$ is homotopy equivalent to $K$.
\end{prop}

The proof of the above (and similar results of \cite{TeVi}) is motivated by our construction of spaces possessing homotopical smallness and homotopical closeness.

\end{document}